\documentclass[11pt, a4paper, reqno]{amsart}
\usepackage{graphicx, xcolor, amsmath, amssymb}

\numberwithin{equation}{section}

\newcommand{\dd}{\mathrm{d}}

\newtheorem{theorem}{Theorem}[section]
\newtheorem{proposition}[theorem]{Proposition}
\newtheorem{definition}[theorem]{Definition}
\newtheorem{lemma}[theorem]{Lemma}

\theoremstyle{remark}
\newtheorem{remark}[theorem]{Remark}
\newtheorem{example}[theorem]{Example}

\title[Propagating fronts for a viscous Hamer-type system]{Propagating fronts for a\\ viscous Hamer-type system} 

\author[G. Cianfarani Carnevale]{Giada Cianfarani Carnevale}
\address[Giada Cianfarani Carnevale]{Dipartimento di Ingegneria e Scienze dell'Informazione e Matematica, Universit\`a degli Studi dell'Aquila (Italy)}
\email{giada.cianfaranicarnevale@graduate.univaq.it}

\author[C. Lattanzio]{Corrado Lattanzio}
\address[Corrado Lattanzio]{Dipartimento di Ingegneria e Scienze dell'Informazione e Matematica, Universit\`a degli Studi dell'Aquila (Italy)}
\email{corrado@univaq.it} 

\author[C. Mascia]{Corrado Mascia}
\address[Corrado Mascia]{Dipartimento di Matematica ``Guido Castelnuovo",
Sapienza Universit\`a di Roma (Italy)}
\email{corrado.mascia@uniroma1.it}

\begin{document}
\baselineskip=16pt

\keywords{}
\subjclass[2010]{}

\maketitle

\begin{abstract}
Motivated by radiation hydrodynamics, we analyse a $2\times2$ system consisting of a one-dimensional viscous conservation law
with strictly convex flux  --the viscous Burgers' equation being a paradigmatic example-- coupled with an elliptic equation,
named {\bf viscous Hamer-type system}.
In the regime of small viscosity and for large shocks, namely when the profile of the corresponding underlying inviscid model undergoes
a discontinuity --usually called {\it sub-shock}-- it is proved the existence of a smooth propagating front, regularising the jump
of the corresponding inviscid equation.
The proof is based on \emph{Geometric Singular Perturbation Theory} (GSPT) as introduced in the pioneering work
of  Fenichel \cite{Fenichel} and subsequently developed by Szmolyan \cite{SZMOLYAN}.
In addition, the case of small shocks and large viscosity is also addressed via a standard bifurcation argument.
\end{abstract}

{\footnotesize
\noindent {\sc Keywords}. parabolic-elliptic system, traveling waves, singular perturbation theory, radiation hydrodynamics.\\
\noindent  {\sc 2020 MSC classification}. 35C07, 34E15, 35B25, 76N30, 35A24.
}







\section{Introduction}

The dynamics of a gas in presence of radiation can be described by the classical compressible Euler equations
with an additional term in the energy balance modelling the radiation effects (see \cite{MM99,Zel02}). 
The extra state variable, describing the intensity of radiation, is often modelled positing that it obeys to an elliptic
equation, leading to the hyperbolic--elliptic system
\begin{equation}\label{generalrad}
	\left\{\begin{aligned}
		\partial_t \rho &+ \partial_x(\rho\,u) = 0\\
		\partial_t(\rho\,u) & + \partial_x(\rho\,u^2 +p) = 0\\
		\partial_t(\rho\,E) & + \partial_x(\rho\,E u +p\,u - \kappa\,\partial_x n) = 0\\
		-\partial_{xx}n & + \tau\left\{n-g(\theta)\right\} = 0
	\end{aligned}\right.
\end{equation}
where $\rho$ denotes density, $u$ velocity, $p$ pressure, $E$ specific total energy and $\theta$ temperature
to be linked by some appropriate constitutive relations.
The classical  Euler system is then coupled through the term $\kappa\,\partial_x n$ to an additional quantity driven by the variable $n$,
describing the average of the photon density by means of the positive parameters $\kappa$ and $\tau$, and the function $g$,
who is usually assumed to have a power-like form.
Such model is obtained in the non relativistic limit of a corresponding hyperbolic--kinetic system where the variable describing
the photon density obeys to a transport equation with interaction kernel given by the Stefan--Boltzmann law
(see \cite{BD04,G-LG05,LMH99} for further details).

A simplified model for radiation dynamics is the $2\times2$ system for the scalars unknowns $u$ and $v$ having the specific form
\begin{equation}\label{Hamoriginal}
	\left\{\begin{aligned}
	\partial_t u+ \partial_x f(u) &= \partial_x v\\
	v - \partial_{xx} v &= \partial_{x} g(u) ,
	\end{aligned}\right.
\end{equation}
for some functions $f$ and $g$ satisfying appropriate assumptions to be specified later.
Here, the inviscid scalar conservation law is augmented with a scalar elliptic equation,
mimicking the coupling present in system \eqref{generalrad}.
For the choices $f(u)=\tfrac{1}{2}u^2$ and $g(u)=u$, system \eqref{Hamoriginal} is known as the {\it Hamer model} for radiating gas,
see \cite{Ham71} (setting $q:=-v$), whose well-posedness is discussed in \cite{LM}.

In the analysis of both systems \eqref{generalrad} and \eqref{Hamoriginal}, an intriguing issue is the study
of existence and stability of shock profiles.
Specifically, the weak dissipation properties --due to the coupling with the elliptic equation-- give rise to the existence
of sub--jumps for sufficiently strong shocks.
Correspondingly, an increasing amount of regularity of the profile emerges as the magnitude of the shock decreases.

Such problem has been addressed in various cases, that is for the ``scalar'' models --a single conservation law with general
flux functions coupled with the elliptic equation-- both for small and regular profiles \cite{ST92}, and for possibly large and
discontinuous ones \cite{Kawashima,LMS,LMS2}. 
The cases of systems has been addressed as well, and in particular for the specific Euler model with radiation effects \eqref{generalrad} in \cite{LCG06}
for weak (regular) profiles, and in \cite{CGLL12,Mas2013} for strong (discontinuous) shocks.
The case of general hyperbolic--elliptic systems for small and large shocks is investigated in \cite{LMS} for linear coupling
and in \cite{LMS2} for nonlinear coupling.

Here, we are interested in exhibiting how the presence of a viscosity term in the simplified Hamer model \eqref{Hamoriginal}
modifies the existence and regularity properties of the shock profiles.
Therefore, in the sequel we consider the following regularized version of system \eqref{Hamoriginal}
\begin{equation}\label{genintro}
\left\{\begin{aligned}
	 \partial_t u + \partial_xf(u) - \epsilon\,\partial_{xx}u  & = \partial_x v \\
	v - \partial_{xx}v  & =  \partial_{x}g(u) ,
\end{aligned}\right.
\end{equation}
where the flux function  $f\in C^2(\mathbb{R})$ is (uniformly) convex
and the coupling function $g\in C^1(\mathbb{R})$ is (strictly) increasing, that is
\begin{equation}\label{hypFG}
	\dd^2 f(u)>0 \quad \textrm{and}\quad \dd g(u)>0
	\qquad\textrm{for any}\; u\in\mathbb{R},
\end{equation}
where $\dd$ denotes the derivative with respect to $u$.
In what follows, we refer to \eqref{genintro} as a {\bf viscous Hamer-type system}.
The corresponding {\it one-field equation} for the unknown $u$ can be obtained by eliminating the variable
$v$ from the second equation and obtaining
\begin{equation*}
	 \partial_t u + \partial_x f(u) - \partial_{xx} \left\{\epsilon\,u+g(u)\right\}
	 	= \partial_{xx} \left\{ \partial_t u + \partial_x f(u) - \epsilon\,\partial_{xx} u\right\}.
\end{equation*}
Incidentally, let us also observe that the method used in \cite{LM} to exhibit well-posedness of the standard Hamer model
passes through a parabolic regularization of the form \eqref{genintro} (see also \cite{CorlRohd12,XuanHuiJin09}).

The main topic of this paper is to investigate the existence of shock profiles for model \eqref{genintro}
in the specific case of end states leading to the presence of a sub--shock in the inviscid model (see \cite{Kawashima,LMS,LMS2}).
\vskip.15cm

For completeness, let us start with a classical definition.

\begin{definition}
Given the states $u_\pm\in\mathbb{R}$, a {\bf propagating front} for the hyperbolic-elliptic system \eqref{genintro}
is a triple $(u,v;c)$ where $(u,v)$ is a travelling wave solution $(u,v)=(u,v)(x-ct)$ satisfying the asymptotic conditions
$(u,v)(\pm\infty) = (u_\pm,0)$ and $c$ is a given constant.

The couple $(u,v)$ is called the {\bf profile} and the parameter $c$ is the {\bf speed}.
\end{definition}

Incidentally, we are denoting with the same symbol both the general solution to \eqref{genintro}, 
and the specific propagating front.
Being this manuscript devoted exclusively to the existence of fronts, we are confident that this will not generate confusion.

The asymptotic states $(u_\pm,v_\pm)$ are forced to be equilibria of \eqref{genintro},
i.e. $u_\pm\in\mathbb{R}$ with $u_+\neq u_-$ and $v_\pm=0$.
Consequently, the speed $c$ is forced to satisfy the {\it Rankine--Hugoniot relation}
\begin{equation}\label{RH}
	c=\frac{f(u_+)-f(u_-)}{u_+ - u_-},
\end{equation}
as it is readily seen by integrating in $\mathbb{R}$ the ordinary differential equation for the profile $(u,v)$.
A standard assumption in the context of scalar conservation laws with convex flux is the {\it Lax condition}
\begin{equation}\label{Lax}
	\dd f(u_+) < c < \dd f(u_-),
\end{equation}
guaranteeing that the jump from $u_-$ is $u_+$ satisfies the entropy condition.
Given the Rankine--Hugoniot condition \eqref{RH} and the convexity assumption on the flux $f$,
inequalities \eqref{Lax} are satisfied if and only if $u_+<u_-$.

Next, we state the two main theorems of our work.
The first one is based on a {\it singular perturbation} approach relative to the (small) parameter $\epsilon$.

\begin{theorem}\label{thm:main1}
Assume \eqref{hypFG}-\eqref{RH}-\eqref{Lax}.
Moreover, let $u_+,u_-$ be such that the inviscid profile undergoes a sub-shock in the u variable.
Then, for $\epsilon >0$ sufficiently small, the parabolic-elliptic system $\eqref{genintro}$ support propagating fronts
with speed $c$ given by the Rankine-Hugoniot relation \eqref{RH}.
\end{theorem}

This result  is obtained in the case of small viscosity coefficient $\epsilon$ and for sufficiently large shocks, namely when
the corresponding non viscous profile experiences a discontinuity, by taking advantage of the Fenichel's \emph{Geometric Singular Perturbation Theory},
in what follows shortened as GSPT (for details, see \cite{Fenichel, SZMOLYAN} and the Appendix at the end of the paper).
The condition which locates the sub--shock in the profiles for the  non viscous models agrees with
the one needed to apply the aforementioned theory, leading to the existence of the connection which smoothes out the discontinuity.

Our second main result, proved in Section \ref{sec:smallshocks}, is based on a direct application of bifurcation theory, where we consider
fixed $O(1)$ viscosity and sufficiently small shocks.

\begin{theorem}\label{thm:main3}
Assume \eqref{hypFG}-\eqref{Lax}. Moreover, let the states $u_+,u_-$ be such that $|u_+-u_-|$ is sufficiently small.
Then for any fixed $\epsilon >0$ the parabolic-elliptic system \eqref{genintro} supports propagating fronts with speed
$c$ given by the Rankine-Hugoniot relation \eqref{RH}.
\end{theorem}

The paper includes a final Appendix \ref{app:GSPT} with, the main general results of GSPT needed in the present analysis.

{\small
\subsection*{Notations}
For readers' convenience, we list (some) symbols used in the text:
\begin{itemize}
\item[] $x$ and $y=x/\epsilon$ denotes the slow and fast variable, respectively;
\item[] $\mathbf{u}=(u,v,w)$ slow/fast variables for $\epsilon\geq 0$;
\item[] $f'$ and $\dot{f}$ denotes the derivatives of the function $f$ with respect to the slow and the fast variable, respectively;
\item[] $\mathcal{S}:=\{\mathbf{u}\,:\,F(\mathbf{u}) = 0\}$ slow manifold with its subsets $\mathcal{S}_\ast$, $\mathcal{S}_\pm$;
\item[] $\mathbf{v}=(u,v,w,\epsilon)$ slow/fast variables;
\item[] $\pi^{\mathcal{S}}$ projection onto the surface $\mathcal{S}$
\item[] $\mathcal{C},\mathcal{C}^s,\mathcal{C}^u$ center/center-stable/center-unstable manifolds for equilibrium point of limiting fast system
\item[] $\mathcal{F}^{s,u} = \lbrace \mathbf{\Psi}^{s,u}(\mathbf{v}) : \mathbf{v} \in \mathcal{C}^{s,u}   \rbrace$ family of stable/unstable manifolds for $\mathcal{C}^{s,u}$
\item[] $W^{s,u}$ stable/unstable manifolds of critical points for reduced slow system
\end{itemize}
However, all the notations will be restated at their first occurrence.
}

\section{General strictly convex case and nonlinear coupling}\label{sec:genconvnon}

To start with, we rephrase the existence of a propagating front in the formalism of singular perturbation theory.

Since $f$ is convex, the function $\dd f$ is strictly increasing and, thanks to \eqref{Lax}, 
there exists a unique point $u_\ast$ such that $\dd f(u_\ast)=c$.
Hence, setting
\begin{equation*}\label{deficit}
	f_c(u):=f(u)-cu.
\end{equation*}
there holds $\dd f_c(u_\ast)=0$.
The function $f_c$ describe the flux in a frame that is co-moving with the propagating front.
The value $f_c(u_+)$ coincides with $f_c(u_-)$ as a consequence of the Rankine--Hugoniot condition \eqref{RH}
and we use the notation $f_c(u_\pm)$ to stress this coincidence.
In addition, the Lax condition becomes
\begin{equation}\label{LaxBis}
	\dd f_c(u_+) < 0 < \dd f_c(u_-).
\end{equation}
For traveling wave solutions, system \eqref{genintro} rewrites as
\begin{equation*}
	\epsilon u''=\left\{f_c(u)-v\right\}',\qquad v''=v-g(u)',
\end{equation*}
where we have used the notation $'= {d}/{dx}$.
Integrating the first equation from $-\infty$ to $x$, and setting $w:=g(u)-v'$,
we obtain the \emph{slow system}
\begin{equation}\label{SLOW SYSTEM}
	\left\{\begin{aligned}
		\epsilon u'&= F(u,v,w):=f_c(u)-f_c(u_\pm)-v\,,\\
	 	v'&= G(u,v,w):=w-g(u)\,,\\
	 	w'&= H(u,v,w):=v\,,\\
	\end{aligned}\right.
\end{equation}
which is a special case belonging to the class of {\it singular perturbed problems}.
The adjective {\it singular} refers to the presence of the parameter $\epsilon$ multiplying first order derivatives of some state
variables and assumed to be small (and positive). 
A significant regime is then obtained as $\epsilon \rightarrow 0$ giving raise to the reduced limit
\begin{equation}\label{SLOW SYSTEM 0}
	\left\{\begin{aligned}	v&=f_c(u)-f_c(u_\pm)\,,\\
	 				v'&=w-g(u)\,,\\
	 				w'&=v\,.	\end{aligned}\right.
\end{equation}
which is called {\it limiting slow system}.
The first equality in \eqref{SLOW SYSTEM 0} can be regarded as an additional constraint along the dynamics,
linking the value of the variable $u$ with the {\it slow variables} $(v,w)$.
As in the present case, the relation $F(u,v,w)=0$ is not invertible with respect to $u$ and the use of appropriate
restrictions have to be considered.

After the rescaling $y:= x/\epsilon$, we get the so-called \emph{fast system}
\begin{equation*}
	\left\{\begin{aligned}
		\dot{u} &= F(u,v,w)\,,\\
	 	\dot{v}&= \epsilon\,G(u,v,w)\,,\\
	 	\dot{w} &= \epsilon\,H(u,v,w)\,.\\
	\end{aligned}\right.
\end{equation*}
with the notation ${}^\cdot = {d}/{dy}= \epsilon\,{d}/{dx}$.
In compact form, setting $\mathbf{u}=(u,v,w)$, the above system can be rewritten as
\begin{equation}\label{FAST SYSTEM}
	\dot{\mathbf{u}}=\mathbf{F}_\epsilon(\mathbf{u})
	:=\bigl(F(\mathbf{u}), \epsilon\,G(\mathbf{u}), \epsilon\,H(\mathbf{u})\bigr).
\end{equation}
As $\epsilon \rightarrow 0^+$ in \eqref{FAST SYSTEM}, we deduce the {\it limiting fast system}
\begin{equation}\label{FAST SYSTEM 0}
	\left\{\begin{aligned}
	\dot u&= f_c(u)-f_c(u_\pm)-v\,,\\
	 \dot v&= 0\,,\\
	 \dot w&= 0\,.\\
	\end{aligned}\right.
\end{equation}
The variable $u$ is said to be the \emph{fast variable} of the system.

From now on, we are confronted with the problem of rigorously establishing the existence of a heteroclinic orbit connecting at $\mp \infty$
the critical points of \eqref{SLOW SYSTEM} --or, equivalently, \eqref{FAST SYSTEM}-- which are 
\begin{equation*}
		\mathbf{u}_\pm =(u_\pm, v_\pm, w_\pm):=\bigl(u_\pm, 0, g(u_\pm)\bigr). 
\end{equation*}
under the assumptions \eqref{hypFG}--\eqref{LaxBis}.
To this aim, we follow the classical singular perturbation approach, consisting in the separate study of dynamical systems
\eqref{SLOW SYSTEM}--\eqref{FAST SYSTEM} in the limiting slow and fast regimes described by \eqref{SLOW SYSTEM 0}--\eqref{FAST SYSTEM 0}, respectively.

Let $\mathcal{S}:=\{\mathbf{u}\in\mathbb{R}^3\,:\,F(\mathbf{u}) = 0\}$ be the \emph{slow manifold} of critical points
for \eqref{FAST SYSTEM 0} and let $\mathcal{S}_\ast$ be the open subset of the slow manifold $\mathcal{S}$ defined by
\begin{equation*}
	\begin{aligned}
	\mathcal{S}_\ast:&=\left\{\mathbf{u}\in\mathcal{S}\,:\, \partial_u F(\mathbf{u})\neq 0 \right\}\\
			&=\mathcal{S} \setminus \left\{\bigl(u_\ast,f_c(u_\ast)-f_c(u_\pm),0\bigr)+t\mathbf{e}_3\,:\,t\in\mathbb{R}\right\},
	\end{aligned}
\end{equation*}
Since the restrictions of the function $\partial_u F$ in the open halflines $(-\infty,u_\ast)$ and $(u_\ast,+\infty)$ are invertible.
system \eqref{SLOW SYSTEM 0} reduces to a two dimensional dynamical system for the slow variables $(v,w)$.
Denoting such inverse functions with the symbols $h_\pm$, we can split the subset $\mathcal{S}_\ast$ as the (disjoint) union
\begin{equation*}
	\mathcal{S}_\ast=\mathcal{S}_-\cup \mathcal{S}_+
\end{equation*}
where
\begin{equation*}
	\begin{aligned}
	\mathcal{S}_-&:=\left\{\mathbf{u}\in\mathcal{S}_\ast\,:\, \partial_u F(\mathbf{u})>0 \right\}\\
			&\,=\left\{(u,v,w)\in\mathbb{R}^3\,:\, v\in\bigl(f_c(u_\ast)-f_c(u_\pm),+\infty\bigr),\, u=h_-(v)\right\},\\
	\mathcal{S}_+&:=\left\{\mathbf{u}\in\mathcal{S}_\ast\,:\, \partial_u F(\mathbf{u})<0 \right\}\\
			&\,=\left\{(u,v,w)\in\mathbb{R}^3\,:\, v\in\bigl(f_c(u_\ast)-f_c(u_\pm),+\infty\bigr),\, u=h_+(v)\right\}.
	\end{aligned}
\end{equation*}
By Lax condition \eqref{LaxBis}, there holds $\mathbf{u}_-\in\mathcal{S}_-$ and $\mathbf{u}_+\in\mathcal{S}_+$.

The analysis for the slow system \eqref{SLOW SYSTEM} exhibits the existence of a discontinuous propagating front for
$|u_+-u_-|$ sufficiently large (in the case of the Hamer model \eqref{Hamoriginal}, the critical threshold is $\sqrt{2}$,
as shown in \cite{Kawashima}).
Specifically, we are interested in showing how the presence of viscosity in the first equation of \eqref{genintro} restores
the regularity of the profile by smoothing out the internal jump of the shock solution.
This task can be accomplished by working on the fast system at $\epsilon=0$ from \cite{Kawashima,LMS,LMS2}
and, then, taking advantage from the results in \cite{SZMOLYAN} to perform a rigorous matching.

\subsection{Splitting consistency} 
The first step consists in the analysis of the linearization at $\mathbf{u}_\pm$ of the fast system \eqref{FAST SYSTEM}
for $\epsilon \geq  0$.
Denoting by $\dd$ the derivative with respect to $\mathbf{u}$, the jacobian of $\mathbf{F}_\epsilon$ is 
\begin{equation*}
	\dd\mathbf{F}_\epsilon (\mathbf{u})
	= \begin{pmatrix} \dd f_c(u) & -1 & 0 \\  -\epsilon\,\dd g(u) & 0 & \epsilon \\  0 & \epsilon & 0 \end{pmatrix}
\end{equation*}
with characteristic polynomial
\begin{equation*}
	\det\left(\dd\mathbf{F}_\epsilon (\mathbf{u})-\lambda\,\mathbb{I}\right)
		=-\lambda^3 + \dd f_c(u)\,\lambda^2 + \epsilon\bigl(\dd g(u)+ \epsilon\bigr)\lambda - \epsilon^2 \dd f_c(u),
\end{equation*}
where $\mathbb{I}$ denotes the identity matrix.
Denoting with $\lambda_1$, $\lambda_2$ and $\lambda_3$ the roots of the polynomial, there hold
\begin{equation*}
	\begin{aligned}
	\textrm{tr}\;\dd\mathbf{F}_\epsilon (\mathbf{u}) 	& = \dd f_c(u)= \lambda_1 + \lambda_2 +\lambda_3\,,\\
	\det \dd\mathbf{F}_\epsilon  (\mathbf{u})		& = -\epsilon^2\dd f_c(u)= \lambda_1 \cdot \lambda_2 \cdot \lambda_3\,. 
	\end{aligned}
\end{equation*}
Next, we analyze the real part of the eigenvalues $\lambda_1$, $\lambda_2$ and $\lambda_3$ of $\dd\mathbf{F}_\epsilon(\textbf{u}_{\pm})$.
Indeed, from the Lax condition \eqref{LaxBis}, it follows
\begin{itemize}
\item[\sl a.] $\textrm{tr}\;\dd\mathbf{F}_\epsilon$ computed at $\mathbf{u}_+$ (respectively $\mathbf{u}_-$) is negative (resp. positive);
\item[\sl b.] $\det\dd\mathbf{F}_\epsilon$ computed at $\mathbf{u}_+$ (resp. $\mathbf{u}_-$)  is positive (resp. negative),\end{itemize}

Since the characteristic polynomial has degree 3, there are two cases: 
\begin{itemize}
\item[\sl i.] either 3 real roots;
\item[\sl ii.] or 1 real root and 2 complex conjugates.
\end{itemize}

Next, we limit the investigation to the case $\mathbf{u}_+$, the other being analogous.

\begin{itemize}
\item[\sl i.] Thanks to {\sl b.}, we have either 3 positive eigenvalues or 2 negative eigenvalues and 1 positive.
The first case is excluded by {\sl a.};
thus, there are 2 negative and 1 positive roots.
\item[\sl ii.] Let us set $\lambda_1\in \mathbb{R}$, $\lambda_2 = \alpha+i\beta $ and $\lambda_3 = \alpha -i\beta$ with $\beta>0$.
Then
\begin{equation*}
	\textrm{tr}\,\dd\mathbf{F}_\epsilon = \lambda_1+2\alpha < 0
	\quad\textrm{and}\quad	
	\det\dd\mathbf{F}_\epsilon =\lambda_1(\alpha^2 + \beta^2) > 0.
\end{equation*}
These conditions imply $\lambda_1>0$ and $\alpha=\textrm{Re}\lambda_{2}=\textrm{Re}\lambda_{3}<0$.
\end{itemize}
Summarizing, at $\mathbf{u}_+$, the matrix $\dd\mathbf{F}_\epsilon$ has a two-dimensional stable manifold and an one-dimensional stable manifold.
At $\textbf{u}_{-}$ the situation is reversed: the dimension of the stable manifold is 1 and that of the unstable manifold is 2.
We infer that the splitting is consistent since the sum of the dimension of the unstable manifold at $\textbf{u}_{-}$ and
the one of the stable manifold at $\textbf{u}_{+}$ is equal to 4 and the state space has dimension equal to 3.

\subsection{Extended fast system}\label{sec:extfs}

Next, we extend the fast system \eqref{FAST SYSTEM} adding a trivial equation for the parameter $\epsilon$, that is
\begin{equation*}\label{layer}
	\dot{\mathbf{v}}=\Phi(\mathbf{v}):=\bigl(F(\mathbf{u}),\epsilon\,G(\mathbf{u}),\epsilon\,H(\mathbf{u}),0\bigr)\,,
\end{equation*}
where $\mathbf{v}=(\mathbf{u},\epsilon)\in\mathbb{R}^3\times(-\epsilon_0,\epsilon_0)$
and $F$, $G$ and $H$ are still defined in \eqref{SLOW SYSTEM}.
Incidentally, let us observe that the above system is equivalent to \eqref{FAST SYSTEM} for $\epsilon>0$ and
to \eqref{FAST SYSTEM 0} for $\epsilon=0$.

The matrix of the linearization of $\Phi$ at $\mathbf{v}$ is
\begin{equation*}
	\dd\Phi(\mathbf{v}) = \begin{pmatrix}
		\dd f_c(u)	 		& -1 		& 0 		& 0 		\\
		-\epsilon\,\dd g(u)	& 0 		& \epsilon 	& w-g(u)	\\
		0				& \epsilon	& 0		& v 		\\
		0				& 0 		& 0		& 0		\end{pmatrix}
\end{equation*} 
with characteristic polynomial at $\mathbf{v}_0:=(\mathbf{u},0)$
\begin{equation*}
	\begin{aligned}
	\det\left(\dd\Phi(\mathbf{v}_0)-\lambda\mathbb{I}\right)
	&=\lambda^3\bigl(\lambda-\dd f_c(u)\bigr).
	\end{aligned}
\end{equation*} 
The polynomial has two roots, $\lambda=0$ and $\lambda=\dd f_c(u)$, which are distinct if $\dd f_c(u)\neq 0$.
The multiplicity 3 of the eigenvalue $0$ corresponds to the dimension of the slow variables $(v,w)$
plus the dimension of the fictitious variable $\epsilon$; the multiplicity 1 of the eigenvalue $\dd f_c(u)$,
coincides with the dimension of the fast variable $u$.
Moreover, since $\lambda =\dd f_c(u) > 0$ on $\mathcal{S}_-$ and $\lambda =\dd f_c(u) <0$ on $\mathcal{S}_+$,
the unstable manifold of $\mathbf{u}_-$ in $\mathcal{S}_-$ and the stable manifold of $\mathbf{u}_-$ in $\mathcal{S}_+$
have both dimension equal to 1.

For readers' convenience, let us recall a basic definition valid for generic compact manifold $M$ and diffeomorphism $f$.

\begin{definition}\label{def:normallyhyperbolic}
An $f$-invariant submanifold $\Lambda$ of $M$ is said to be {\bf normally hyperbolic} if the restriction $T_\Lambda M$ to $\Lambda$
of the tangent bundle $TM$ admits a splitting into a sum of three $\dd f$-invariant sub-bundles,
\begin{align*}
	T_\Lambda M = T\Lambda + E^s +E^u.
\end{align*}
where $T\Lambda$ is the tangent bundle of $\Lambda$ and $E^{s,u}$ denotes the stable/unstable bundle, respectively.
\end{definition}

Here, we consider compact subsets $\mathcal{K}\Subset \mathcal{S}_0$ that are normally hyperbolic invariant manifolds
of the layer problem \eqref{FAST SYSTEM 0}. 

Since $\mathbf{F}_0$ restricted to $\mathcal{S}_\pm$ is identically zero, there holds
\begin{equation*}
	T_{\mathbf{u}}\mathcal{S}_\pm=\ker \dd\mathbf{F}_0.
\end{equation*}
The subspace $T_{\mathbf{u}}\mathcal{S}_\pm$ is invariant under the action of $\dd\mathbf{F}_0$ and, therefore,
the linear map induced by the linearization
\begin{equation*}
	Q\mathbf{F}_0(\mathbf{u})\,:\,T_\mathbf{u} \mathbb{R}^3 / T_\mathbf{u} \mathcal{S}_\pm
	\rightarrow T_\mathbf{u} \mathbb{R}^3 / T_\mathbf{u} \mathcal{S}_\pm
\end{equation*}
is well defined.
 
The eigenvalues of $\dd\mathbf{F}_0$  are exactly the nontrivial ones cited before, since
\begin{equation*}
	\dd\mathbf{F}_0(\mathbf{u}) = \begin{pmatrix} (\dd f_c\circ h_\pm)(v) & -1 & 0 \\ 0 & 0 & 0 \\  0 & 0 & 0  \end{pmatrix}.
\end{equation*}
Moreover, a generic element $\mathbf{v}\in T_{\mathbf{u}}\mathcal{S}_\pm$ is $\mathbf{v} = \alpha\,\bar{\mathbf{u}}+\beta\,\mathbf{e}_3$
for $(\alpha,\beta)\in\mathbb{R}^2$ where $\bar{\mathbf{u}}:=\left(1,(\dd f_c\circ h_\pm)(v),0\right)$.

Next, following \cite{SZMOLYAN}, we define the {\it projection map} $\pi^\mathcal{S}$ by the splitting 
\begin{equation*}
	T\,\mathbb{R}^3 = T\mathcal{S}_\pm \oplus N,
\end{equation*}
where $N$ is the complement of $T\mathcal{S}_\pm$ invariant under $\dd\mathbf{F}_0$.
The matrix $\dd\mathbf{F}_0$ maps the vector $\mathbf{v}=(a,b,c)\in\mathbb{R}^3$ into
\begin{equation*}
	\begin{pmatrix} (\dd f_c\circ h_\pm)(v) & -1 & 0 \\ 0 & 0 & 0 \\  0 & 0 & 0  \end{pmatrix}
	\begin{pmatrix} a \\ b \\ c \end{pmatrix}
	= \begin{pmatrix} (\dd f_c\circ h_\pm)(v) a - b \\ 0 \\ 0 \end{pmatrix}.
\end{equation*}
Therefore, for a nonzero vector $(a,b,c)$ to be invariant under the action of $J\mathbf{F}_0$ it is required $b=c=0$.
Since $(\dd f_c\circ h_\pm)(v)\neq 0$, the subspace $N$ is thus spanned by $\mathbf{e}_1:=(1,0,0)$.
Summarizing, the splitting is
\begin{equation*}
	T\,\mathbb{R}^3 \bigr|_{\mathcal{S}_\pm} = \langle \bar{\mathbf{u}}, \mathbf{e}_3 \rangle \oplus \langle \mathbf{e}_1 \rangle .
\end{equation*}
Let us denote with $\mathcal{C}$, $\mathcal{C}^s$ and $\mathcal{C}^u$ and  the center, center-stable, center-unstable manifolds
generated by the extended fast system $\Phi$.
They are locally invariant manifolds containing $\mathcal{K} \times \{ 0 \}$ and tangent to the corresponding center, center-stable,
center-unstable eigenspaces of the linearization $\dd\Phi$ at $\mathbf{v}=(\mathbf{u},0) \in \mathcal{K} \times \{0\}$,
denoted here by $E_{\mathbf{u}}^c$, $E_{\mathbf{u}}^s \oplus E_{\mathbf{u}}^c$ and $E_{\mathbf{u}}^u \oplus E_{\mathbf{u}}^c$,

In the case under scrutiny, the dimension $k^c$ of the center manifold for the system \eqref{FAST SYSTEM} is 0.
We recall that $n=2$ is the dimension of the slow variable $(v,w)$ and $k=k^u+k^s=1$ equals the dimension of the fast variable.
The following statements are true:
\begin{itemize}
\item[1.] $\dim\,\mathcal{C}=n+1= 3$;
\item[2.] for the surface $\mathcal{S}_+$, $k^s=1$ and $k^u=0$, so that $\dim \mathcal{C}^s=n+1+k^s=4$
	and $\dim\,\mathcal{C}^u=n+1+k^u=3$;
\item[3.] for the surface $\mathcal{S}_-$, $k^s=0$ and $k^u=1$ so that $\dim\,\mathcal{C}^s=n+1+k^s=3$
	and $\dim\,\mathcal{C}^u=n+1+k^u=4$.
\end{itemize}

Given constants $p, q\in[1,+\infty)$, let $\mathcal{C}^s$ be a center-stable manifold for $\Phi$
in a neighborhood of $\mathcal{K}\times\{0\}$.
At this point, we need to introduce also the notion of \emph{family of stable (resp. unstable) manifolds}
for $\mathcal{C}^s$ (resp. $\mathcal{C}^u$).

\begin{definition}
A family $\mathcal{F}^s=\left\{\Psi^s(\mathbf{v})\,:\, \mathbf{v} \in \mathcal{C}^s \right\}$ is a {\bf $C^{q}$-family
of $C^{p}$-stable manifold for $\mathcal{C}^s$ near $\mathcal{K}$} if
\begin{itemize}
\item[i.] $\Psi^s(\mathbf{v})$ is $C^{p}$-manifold for each $\mathbf{v} \in \mathcal{C}^s$;
\item[ii.] $\mathbf{v} \in \Psi^s(\mathbf{v})$ for each $\mathbf{v} \in \mathcal{C}^s$;
\item[iii.] $\Psi^s(\mathbf{v}_1)$ and $\Psi^s(\mathbf{v}_2)$ are disjoint if $\mathbf{v}_1\neq \mathbf{v}_2$
	for each $\mathbf{v}_1,\mathbf{v}_2 \in \mathcal{C}^s$;
\item[iv.] $\Psi^s(\mathbf{u},0)$ is tangent to $E_{\mathbf{u}}^s$ at $(\mathbf{u},0)$ for each $\mathbf{u} \in \mathcal{S}$;
\item[v.] the set $\lbrace \Psi^s(\mathbf{v})\,:\,\mathbf{v} \in \mathcal{C}^s \rbrace$ is a positively invariant $C^{q}$-family of manifolds
with respect to the flow $\Phi$.
\end{itemize}
\end{definition}

The case of the unstable manifolds $\mathcal{F}^u$ is defined in an analogous way.
These families provide a $\emph{foliation}$ of $\mathcal{C}^s$ and $\mathcal{C}^u$, i.e.
\begin{equation*}
	\mathcal{C}^s = \left\{\Psi^s(\mathbf{v}): \mathbf{v} \in \mathcal{C} \right\}
	\quad\textrm{and}\quad
	\mathcal{C}^u = \left\{\Psi^u(\mathbf{v}): \mathbf{v} \in \mathcal{C} \right\}.
\end{equation*}

\subsection{Construction of the profile for the reduced system}\label{sec:red0}
We want to apply Lemma \ref{lemma3.1} in  {Appendix A} to system \eqref{FAST SYSTEM} to describe
the flow induced on $\mathcal{S}$ when it is given by the graph of a function, i.e. in $\mathcal{S}_\pm$.

Adapting Lemma 5.4 in \cite{Fenichel} to \eqref{FAST SYSTEM}, we infer that
the projections $\pi^\pm$ on the surfaces $\mathcal{S}_{\pm}$ are given by the multiplication against the matrices
\begin{equation*}
	\mathbb{A}^\pm  = \begin{pmatrix}
		0	& \dd h_\pm(v)	& 0	\\
		0	& 1 			& 0	\\
		0	& 0			& 1
				\end{pmatrix},
\end{equation*}
so that the \emph{reduced systems} are
\begin{equation*}
	\mathbf{u}' = \mathbb{A}^\pm \begin{pmatrix} 0 \\ G \\ H \end{pmatrix}
	= \begin{pmatrix} \dd h_\pm(v) \,G(h_\pm(v), v,w) \\ G(h_\pm(v), v,w) \\ H(h_\pm(v), v,w) \end{pmatrix}
\end{equation*}
In components, for any $(u,v,w) \in \mathcal{S}_{\pm}$, the above systems become
\begin{equation}\label{reduced}
	\left\{\begin{aligned}
		u' &= \frac{w-(g\circ h_{\pm})(v)}{(\dd f_c\circ h_{\pm})(v)} \\
	 	v'&= w-(g\circ h_{\pm})(v) \\
	 	w'&= v. \\
		\end{aligned}\right.
\end{equation}
In what follows, we refer to the vector field as $\mathbf{F}_R$. 

\begin{remark}
In the case under investigation, the computations yielding to \eqref{reduced} reduces to differentiate \eqref{SLOW SYSTEM 0}$_1$
along the profile lying on $\mathcal{S}_\pm$. Hence, for $u= h_{\pm}(v)$, and using \eqref{SLOW SYSTEM 0}$_2$ we end up with
\[
	u' = \frac{w-(g \circ h_{\pm})(v)}{( \dd f_c \circ h_{\pm})(v)}
\]
and thus  \eqref{reduced}.
On the other hand, Lemma \ref{lemma3.1} in  {Appendix A} refers to a general framework, where the vector field in \eqref{FAST SYSTEM}
merely  depends  in a regular way on $\epsilon$.
\end{remark}

Upon observation, it is readily seen that the critical points of  the \emph{limiting slow system} \eqref{SLOW SYSTEM 0},
or equivalently \eqref{reduced}, where we have desingularized the problem, are
\begin{equation*}
	\textbf{u}_\pm := \bigl(u_\pm,0,g(u_\pm)\bigr)\in \mathcal{S}_\pm,
\end{equation*}
and hence they concide with the ones of  the perturbed slow systems \eqref{SLOW SYSTEM}.
Theorem \ref{teo3.3} in Appendix guarantees the existence for $\epsilon > 0$ sufficiently small of a sequence of critical points $\mathbf{u}_\epsilon$
for the perturbed vector field in \eqref{SLOW SYSTEM}, smoothly depending on $\epsilon$, and which reduces for $\epsilon =0$ to the critical points
of  \eqref{SLOW SYSTEM 0}.
In our case, again in view of the independence from $\epsilon$ of the vector field  in \eqref{SLOW SYSTEM},
we have $\mathbf{u}_\pm = \mathbf{u} = \mathbf{u}_\epsilon$, without invoking Theorem   \ref{teo3.3}.
However, in order to define the local stable/unstable manifolds of $\mathbf{u}_\epsilon$ needed for the connection,
we shall refer to the general framework of Theorem \ref{s-u-mxr}
and therefore one has to check that $\lambda=1$ is not an eigenvalue of the Jacobian $\dd \mathbf{F}_R(\textbf{u}_\pm)$.
To this end, we compute the Jacobian of the vector fields in \eqref{reduced}
\begin{equation*}
	\dd \mathbf{F}_R(\textbf{u}_\pm)
	=\begin{pmatrix} 	0 & -\dd g/\dd f_c^2 & 1/\dd f_c	\\ 
					0 & -\dd g/\dd f_c^2 & 1 			\\ 
					0 & 1 			 & 0			\end{pmatrix}
\end{equation*}
where $\dd f_c$ and $\dd g$ are computed at $u_\pm$.
The corresponding characteristic polynomials $p_\pm$ are
\begin{equation*}
	p_\pm(\lambda):=-\lambda\, q_\pm(\lambda)
	\quad\textrm{with}\quad
	q_\pm(\lambda)=\lambda^2 + \frac{\dd g(u_{\pm})}{\dd f_c(u_{\pm})}\,\lambda -1.
\end{equation*}
These have three roots: the trivial one, given by $\lambda_\pm^0=0$, and the roots of $q_\pm$,
which can be expressed as
\begin{equation*}
	\begin{aligned}
		\lambda_\pm^{1} & =  -\frac{\sqrt{\dd g^2+4\,\dd f_c^2}+\left(\textrm{sgn}\,\dd f_c\right) \dd g}{2|\dd f_c|}, \\
		\lambda_\pm^{2} & =  \frac{\sqrt{\dd g^2+4\,\dd f_c^2}-\left(\textrm{sgn}\,\dd f_c\right) \dd g}{2|\dd f_c|},
	\end{aligned}
\end{equation*}
where $\dd f_c$ and $\dd g$ are computed at $u_\pm$.
Recalling the Lax condition \eqref{LaxBis} and the monotonicity assumption \eqref{hypFG} on function $g$,
a direct computation shows that
\begin{equation*}
	\lambda_-^1 < -1 < \lambda_+^1 < 0 < \lambda_-^2 < 1 < \lambda_+^2,
\end{equation*}
so that, in particular, $\lambda_\pm^{2}\neq 1$.

Therefore we can apply Theorem \ref{s-u-mxr} to characterize the local stable-unstable manifolds of the perturbed critical points.
These manifolds are of dimension 1 for both $\textbf{u}_\pm$ and thus the critical points for the perturbed system are saddles in the slow directions.

Now we pass to the study of the reduced system \eqref{reduced} to briefly recast the results by Lattanzio et al. \cite{LMS, LMS2}
in the present framework (see also \cite{Kawashima}).
To this end, we focus on the case of sufficiently large shocks so that the profile exhibits a sub--shock for the variable $u$, 
case for which GSPT can be directly applied, because the two branches of the profile belong to $\mathcal{S}_-$ and $\mathcal{S}_+$, respectively.

In the case of the \emph{singularly perturbed systems} \eqref{SLOW SYSTEM},
we prove the existence of the heteroclinic orbit as a perturbation of the so-called \emph{singular heteroclinic orbit}.
The latter consists of orbits of the reduced systems lying  on the two surfaces $\mathcal{S}_\pm$ of the critical manifold $\mathcal{S}$,
linked by a heteroclinic connection of the \emph{layer system} \eqref{FAST SYSTEM 0}.

For completeness, we adapt here the results proved in  \cite{LMS} for system \eqref{genintro} with $g(u)=u$
and in \cite{LMS2} for general increasing functions $g$.
For readers' convenience, we observe that the variable $(u,v,w)$ of the present paper corresponds to $(u,-q,z)$ of \cite{LMS,LMS2}.

Let $(v_\pm,w_\pm)$ be the \underline{maximal} solutions to the problems 
\begin{equation*}\label{wpm}
	v'= w-(g\circ h_{\pm})(v)\,,\qquad w'= v\,,
\end{equation*}
with the asymptotic conditions $(v_\pm,w_\pm)(\pm\infty)=(0,g(u_\pm))$,
describing the dynamics of the reduced system \eqref{reduced} along the surfaces $\mathcal{S}_\pm$. 

Functions $v_-$, $w_\pm$ are decreasing and $v_+$ is increasing.
In addition, there exist $x_\pm\in\mathbb{R}$ such that
\begin{equation*}
	w_\pm(x_\pm) -v'_\pm(x_\pm) = g(u_\ast)  \quad\textrm{and}\quad v_\pm(x_\pm) = -f_c(u_\pm)
\end{equation*}
In particular, the solution $(v_-,w_-)$ is defined in $(-\infty,x_\ast]$ and the solution $(v_+,w_+)$ is defined in $[x_\ast,+\infty)$.

As a consequence of the monotonicity of functions $v_\pm$, there holds
\begin{equation*}
	w_-(x_\ast)\leq g(u_\ast) \leq w_+(x_\ast).
\end{equation*}
The complete heteroclinic orbit is built by matching two branches of the maximal solutions $(v_\pm,w_\pm)$
at a given point $x_\ast$ which can be arbitrarily chosen taking advantage of translation invariance.

\begin{proposition}\label{teo3.5}
Assume hypotheses \eqref{hypFG}.
Then, there exists a function $(v,w)$ with $v\in C^0(\mathbb{R})\cap C^1\left(\mathbb{R}\setminus\{x_\ast\}\right)$
and $w\in C^1(\mathbb{R})\cap C^2\left(\mathbb{R}\setminus\{x_\ast\}\right)$ solving the first order system (not in normal form)
\begin{equation*}
	(f_c\circ g^{-1})(w-v')=v\,,\qquad w'=v\,,
\end{equation*}
satisfying the asymptotic conditions $(v,w)(\pm \infty) = \left(0,g(u_{\pm})\right)$.
Moreover, the function $v'=w''$ has at most a jump discontinuity at $x_\ast$.
\end{proposition}

The existence theorem of propagating fronts for system \eqref{genintro} with $\epsilon=0$
is a simple consequence of the above construction.

\begin{theorem}\label{teo:existsinglim}
Assume hypotheses \eqref{hypFG}.
For any $u_\pm$ satisfying the Lax condition \eqref{Lax}, there exists a propagating front for system \eqref{genintro} with $\epsilon = 0$.
The profile $(u,v)$ is unique up to translation and the speed $c$ is given by the Rankine--Hugoniot condition \eqref{RH}.

Moreover, the component $u$ belongs to $C^1(\mathbb{R}\setminus\{x_\ast\})$, eventually with a jump discontinuity at $x_\ast$.
\end{theorem}

For the complete proofs of Proposition \ref{teo3.5} and Theorem \ref{teo:existsinglim}, we refer to \cite{LMS, LMS2}.
Here, we limit to sketch the basic steps in the matching procedure with the aim of providing support to the subsequent parts.

We define a profile $w$ by matching together the two maximal solutions
\begin{equation*}
	w(x) = \begin{cases} w_-(x) & x< x_\ast 		\\
					w_+(x) & x\geq x_\ast,	\end{cases}
\end{equation*}
relying on the monotonicity of maximal solutions $w_\pm$, such that $w(x_\ast)=w_\ast$.
The latter equality can be written also as
\begin{equation*}
	g(u_r)+v'_r =g(u_\ell)+v'_\ell.
\end{equation*} 
Analogously, the continuity of $w'$ at $x_\ast$ gives
\begin{equation*}
	w'_+(x_\ast)=v_r=v_\ell=w'_-(x_\ast).
\end{equation*}
Since we are investigating a discontinuous profile $u$, we know that
\begin{equation}\label{vastnegative}
	v_\ast:=v_\ell=v_r<0,
\end{equation}
and 
\begin{equation}\label{eq:ulurdiff}
	u_r 
		= g^{-1}\left(w_\ast - v'_+(x_\ast)\right) < u_\ast < g^{-1}\left(w_\ast - v'_-(x_\ast)\right)
		= u_\ell.
\end{equation}
Therefore, there exists $\delta > 0$ sufficiently small such that, along the constructed profile, we have
\begin{equation*}
	h_+(v) \in [u_+,u_\ast-\delta], \quad
	h_-(v) \in [u_\ast + \delta, u_-],
\end{equation*}
and GSPT is applicable.
\begin{remark}
The present analysis is valid also for the particular case of the viscous Hamer model with linear coupling \eqref{Hamoriginal},
namely for $f(u)= u^2/2$, $g(u)=u$.
In that case,  the slow manifolds are given by: 
\begin{equation*}
	\mathcal{S}_\pm = \lbrace (u,v,w) \in \mathbb{R}^3 \ / \ u = \mp \sqrt{1+2v} \rbrace.
\end{equation*}
If we choose the asymptotic states as $u_\pm =(\mp 1,0, \mp 1)$, we obtain $u^*=0$ and $c=0$.
In \cite{Kawashima, LMS}, the authors showed that the profile increases its regularity  as the strength of the shock decreases.
Moreover, for that specific model, from \cite{Kawashima} we know that 
, $u$ has precisely an ammissible jump discontinuity at a single  point  if $|u_+-u_-| > \sqrt 2$, which is the case we are referring to here. 
Therefore we get only the $C^1$-continuity of $w$. The maximal solution $w_+$ decreases from zero toward -1 at $+\infty$
and $w_-$ decreases from 1 at $-\infty$ toward 0. Due to the particular choice of the flux $f$, and thus of the inverse function $h_\pm(v)$,
we get $u_\ell = -u_r$. We shall  give a partial answer in the case of small schocks  in Section \ref{sec:smallshocks},
when we shall discuss the existence of the profile for sufficiently small shocks and arbitrary diffusion.
\end{remark}

\subsection{Construction and persistence of the orbit}\label{subsec:singular}

In the previous section we have constructed the first two pieces of the singular heteroclinic orbit,
namely orbits of the reduced systems $\mathbf{F}_R$ lying on $\mathcal{S}_{\pm}$.
Here, we aim to substitute the istantaneous jump for $u$ among these branches with an appropriate smoothed
version described by the limiting equation \eqref{SLOW SYSTEM 0} of the original differential equations
\eqref{SLOW SYSTEM} in the regime $\epsilon\to 0^+$.

The limiting fast system \eqref{FAST SYSTEM 0} has critical points given by $(h_{\pm}(v),v,w)$ for arbitrary $v, w$.
Therefore, we make use of it to connect the two manifolds $\mathcal{S}_{\pm}$ at points ${\mathbf u}_\ell=(u_\ell,v_\ast,w_\ast)\in\mathcal{S}_-$
and ${\mathbf u}_r = (u_r,v_\ast,w_\ast)\in\mathcal{S}_+$, where $u_{\ell}= h_-(v_\ast)$ and $u_r=h_+(v_\ast)$
The dynamics are completely described by the reduced equation
\begin{equation}\label{eq:ODE1non}
	\dot{u} = F(u,v_\ast,w) = f_c(u)-f_c(u_\pm)- v_\ast.
\end{equation}
Thus, for fixed $v_\ast>f_c(u_\ast)-f_c(u_-)$, it is easily checked that
\begin{equation*}
	\begin{aligned}
	&F(u,v_\ast,w)= 0\quad\textrm{if and only if}\quad u=u_{\ell}\quad\textrm{or}\quad u=u_{r}\,, \\
	&F(u,v_\ast,w)< 0\quad\textrm{for}\quad u \in (u_r,u_\ell)\,,
	\end{aligned}
\end{equation*}	
with the opposite sign for $u\in (- \infty,u_r)\cup(u_\ell,\infty)$.
We infer the existence of a global and decreasing solution $u_0$ to \eqref{eq:ODE1non} verifying the asymptotic
conditions $u_0(-\infty)=u_{\ell}$ and $u_0(+\infty)=u_{r}$.

\begin{example}[Hamer model]\label{rem:Hamer_u0}
Let $u_\mp =\pm 1$.
As a consequence of the Rankine--Hugoniot \eqref{RH}, the speed $c$ is zero and
\begin{equation*}
	f_0(u)=f(u)-0\cdot u=\tfrac12 u^2.
\end{equation*}
Moreover, since $|u_- - u_+| = 2 >\sqrt 2$, the inviscid profile undergoes a Lax sub--shock  between $u_\ell$ and $u_r = - u_\ell <0$. 
For any fixed value $v_\ast\in\left(-\tfrac12 ,0\right)$, we can explicitly solve  the equation \eqref{eq:ODE1non} for $u$.
Indeed, setting $v_\ast:=-\tfrac12(u_-^2-1)$ for some $u_\ell\in(0,1)$, there holds
\begin{equation*}
	\dot{u}=   \tfrac{1}{2}(u^2-1)- v_\ast = \tfrac{1}{2}(u^2-u_\ell^2).
\end{equation*}
Separating the variables and integrating by parts, we obtain
\begin{equation*}
	u_0(x)= - u_\ell\cdot\frac{e^{u_\ell x}-1}{e^{u_\ell x}+1}=-u_\ell \tanh\Bigl(\tfrac12u_\ell x\Bigr),
\end{equation*}
which satisfies the required asymptotic conditions $u_0(\pm\infty)=\mp u_\ell$.
\end{example}

Let $a, b\in\mathbb{R}$ with $a\leq b$.
Then, given $\delta\in\bigl(0,f_c(u_\pm)\bigr)$, the compact normally hyperbolic invariant
manifolds $\mathcal{K}_\pm$  have the form 
\begin{equation*}
	\mathcal{K}_\pm = \left\{ \mathbf{u}\in \mathcal{S}_\pm\,:\,(v,w)\in \bigl[\delta - f_c(u_\pm),0\bigr]\times[a,b]\right\}.
\end{equation*} 
It is clearly seen that the two parameter sets coincide.

Let $\mathcal{W}_{-,0}^u$ be the unstable manifold of $\textbf{u}_-$ and $\mathcal{W}_{+,0}^s$ the stable manifold of $\textbf{u}_+$
with respect to the reduced slow system \eqref{reduced}.
Then, the {\it singular unstable manifold} of $\textbf{u}_-$ and {\it singular stable manifold} of $\textbf{u}_+$ are given by
\begin{equation*}
	N_{-,0}^u = \bigcup_{\mathbf{u} \in \mathcal{W}_{-,0}^u} \mathbf{\Psi}^u(\mathbf{u})
	\quad\textrm{and}\quad
	N_{+,0}^s = \bigcup_{\mathbf{u} \in \mathcal{W}_{+,0}^s}  \mathbf{\Psi}^s(\mathbf{u}),
\end{equation*}
where $\mathbf{\Psi}^{u,s}(\mathbf{u})$ are the unstable/stable fibers of the limiting fast system \eqref{FAST SYSTEM 0} based at $\mathbf{u}$,
which are tangent to the eigenvectors associated to the non-trivial eigenvalue of the Jacobian $\dd\mathbf{F}_0(\textbf{u})$
with an appropriate sign choice. 
The construction of the singular heteroclinic orbit is complete if and only if we show that the intersection of the manifolds $N_{-,0}^u$ and $N_{+,0}^s$
is transversal along the solution of \eqref{eq:ODE1non}, namely upon construction of the smooth connection between $u_\ell$ and $u_r$ recalled above.

For $\mathbf{u}\in\mathcal{S}_\ast$, the non trivial eigenvalue of the Jacobian $\dd \mathbf{F}_0(\textbf{u})$ is
$\lambda =\dd f_c(u)$ which is negative in $\mathcal{S}_+$ and positive in $\mathcal{S}_-$.
At $\mathbf{u} \in \mathcal{S}_-$, the eigenvector relative to the eigenvalue $\dd f_c(u)$ is readily seen to be $\mathbf{e}_1=(1,0,0)$,
so that the unstable fiber $\Psi^u(\mathbf{u})$ is explicitly given by $\mathbf{u}+t\mathbf{e}_1$ for $t\in\mathbb{R}$.
The same computation is valid for $\mathbf{u} \in \mathcal{S}_+$ so that the elements of $\Psi^s(\mathbf{u})$
have the form $\mathbf{u}+t\mathbf{e}_1$ for some $t\in\mathbb{R}$.

\begin{example}[Hamer model]
For the viscous Hamer model, the sets $\mathcal{K}_{\pm}$ have the form 
\begin{equation*}
	\mathcal{K}_\pm = \left\{ \mathbf{u}\in \mathcal{S}_\pm\,:\,(v,w)\in \left[\delta - \tfrac{1}{2},0 \right]\times[a,b]\right\}.
\end{equation*}
for $\delta \in \left(0, \frac{1}{2} \right)$. The unstable/stable fibers $\Psi^{u,s}(\mathbf{u})$ at $u_\pm = (\mp 1, 0, \mp 1)$ are also tangent to $\mathbf{e}_1$.
\end{example}

In view of the previous discussion, for $v_\ast>f_c(u_\ast)-f_c(u_-)$, there exists a connection $u_0$, solution to \eqref{eq:ODE1non}, between the  points
\begin{equation*}
	\mathbf{u}_\ell = (u_\ell, v_\ast,w_\ast) \in \mathcal{W}_{-,0}^u\subset \mathcal{K}_-
	\quad\textrm{and}\quad
	\mathbf{u}_r =  (u_r, v_\ast,w_\ast) \in \mathcal{W}_{+,0}^s\subset\mathcal{K}_+,
\end{equation*}
where $u_{l}= h_-(v_\ast)$ and $u_r= h_+(v_\ast)$.
Such connections are uniquely determined in view of the monotonicity of the component $w$
along the stable/unstable manifolds $\mathcal{W}_{+,0}^s$ and $\mathcal{W}_{-,0}^u$.

The subsequent Theorem --an ultra-simplified version of \cite[Theorem 4.1]{SZMOLYAN} which for completeness we stated
in the Appendix (see Theorem \ref{teo3.12})-- gives conditions to prove transversality in the regime $\epsilon = 0$.
\begin{theorem}\label{thm:trans}
Let $\phi(\mathcal{W}_{-,0}^u)$ and $\phi(\mathcal{W}_{+,0}^s)$ denote the $(v,w)$-coordinates of the manifolds
$\mathcal{W}_{-,0}^u\subset\mathcal{S}_-$ and $\mathcal{W}_{+,0}^s\subset\mathcal{S}_+$, respectively.
Then, the manifolds $N_{-,0}^u$ and $N_{+,0}^s$ intersect transversally at the points of the heteroclinic orbit if and only if
\begin{equation}\label{nointer}
	T_{(v_\ast,w_\ast)}\phi(\mathcal{W}_{-,0}^u) \cap T_{(v_\ast,w_\ast)}\phi(\mathcal{W}_{+,0}^s) = \{ 0 \}.
\end{equation}
\end{theorem}
In the sequel, we check the validity of \eqref{nointer}.

\begin{remark}
The intersection of $N_{-,0}^u$ and $N_{+,0}^s$ at a point of the form $\mathbf{u}=(u_0(\bar{x}),v_\ast,w_\ast)$ (for some generic
$\bar x\in\mathbb{R}$) is transversal if and only if $T_{\mathbf{u}} N_{-,0}^u + T_{\mathbf{u}} N_{+,0}^s = \mathbb{R}^3$,
that is if and only if
\begin{equation*}
	\dim(T_{\mathbf{u}} N_{-,0}^u + T_{\mathbf{u}} N_{+,0}^s)=3.
\end{equation*}
Following the notation in \cite{SZMOLYAN}, the equality 
\begin{equation*}
	\begin{aligned}
	\dim(T_{\mathbf{u}} N_{-,0}^u + T_{\mathbf{u}} N_{+,0}^s) & = \dim(T_{\mathbf{u}} N_{-,0}^u)
	 + \dim(T_{\mathbf{u}} N_{+,0}^s)\\
	& \ - \dim(T_{\mathbf{u}} N_{-,0}^u \cap T_{\mathbf{u}} N_{+,0}^s)
	\end{aligned}
\end{equation*}
implies that the intersection in \eqref{nointer} is transversal if and only if
\begin{equation*}
	\dim(T_{\mathbf{u}} N_-^u \cap T_{\mathbf{u}} N_+^s) = d,
\end{equation*}
where $d$ in our case is given by
\begin{equation*}
	d= j_-^u+k_-^u+ j_+^s+k_+^s - n-k = 1 + 1 + 1 + 1 - 2 - 1 = 1,
\end{equation*}
with $j_-^u$, $k_-^s$, $j_+^s$ and $k_+^s$ describing the dimensions of the stable and unstable reduced manifolds
of the singular points $\mathbf{u}_\pm$.
In the general setting, Theorem \ref{teo3.12} guarantees the needed transversal intersection along  points of the heteroclinic orbit if and only if there 
exist exactly $d-1$ linearly independent solutions $\xi \in T_{(v_\ast,w_\ast)}\phi(W_{-,0}^u) \bigcap T_{(v_\ast,w_\ast)}\phi(W_{+,0}^s)$ of the equation
\begin{equation*}
	(M,\xi)= 0,
\end{equation*}
where $M \in \mathbb{R}^2$ is defined by
\begin{equation*}
	M: = \int_{\mathbb{R}} \psi(\xi)\,\partial_{v,w} \bigl\{f_c(u)-f_c(u_-)+v\bigr\}\,d\xi
\end{equation*}
and $\psi$ is defined as the unique bounded solution of 
\begin{equation*}\label{eq:odeadded}
	\psi' = -\dd f_c(u_0(x))\,\psi.
\end{equation*}
Here $d-1=0$ and  $\partial_{v,w} (f_c(u)-f_c(u_-)+v)= (1,0)$.
Moreover
\begin{equation*}
	\psi(x) = \psi(0)\exp\left\{-\int_0^{x} \dd f_c(u_0)\,d\xi\right\}.
\end{equation*}
Since $u_r < u_0(x) <u_\ell$ and $\dd f_c$ is monotone, then 
\begin{equation*}
	\dd f_c(u_r)< \dd f_c (u_0(x)) < \dd f_c(u_\ell),
\end{equation*}
which implies 
\begin{equation*}
	\psi(0)e^{-\dd f_c(u_\ell) x} < \psi(0)e^{-\int_0^{x} \dd f_c(u_0) ds} < \psi(0)e^{-\dd f_c(u_r) x}.
\end{equation*}
In view of the conditions  $\dd f_c(u_r)< 0<\dd f_c(u_\ell)$, we can conclude that  the only $\psi$ which is globally
bounded in $\mathbb{R}$ is the trivial one, that is $\psi\equiv 0$ and as a consequence the vector $M$ is identically 0.
In other words, in the present case, the requirements of Theorem  \ref{teo3.12} reduce to \eqref{nointer}.
\end{remark}

Analyzing the intersection in \eqref{nointer}, this condition is equivalent to inquire the existence of a couple 
$(v_0, w_0) \in R:= [\delta- f_c(u_\pm), 0]  \times  [a,b] \subset \mathbb{R}^2$
such that the points $\left(h_-(v),v,w\right) \in \mathcal{W}_{-,0}^u$ and $\left(h_+(v),v,w\right) \in \mathcal{W}_{+,0}^s$
are connected by a heteroclinic orbit of the layer problem \eqref{FAST SYSTEM 0}, i.e. a solution $u_0$
to \eqref{eq:ODE1non} for some constant $w\in\mathbb{R}$.
By construction, both the maximal solution $w_-$ (describing the unstable manifold $\mathcal{W}_{-,0}^{u}$)
and the maximal solution $w_+$ (describing the stable manifold $\mathcal{W}_{+,0}^{s}$) are monotone decreasing.
Hence, since $w_-$ decreases from $g(u_-)$ at $-\infty$ toward $w_\ast$ and $w_+$ decreases from $w_\ast$ toward $g(u_+)$ at $+\infty$,
we conclude that the only point in $R$ verifying the above conditions is
\begin{equation*}
	(v_0,w_0)=(v_r,w_\ast)=(v_\ell,w_\ast)=(v_\ast, w_\ast).
\end{equation*}
In other words, in order to construct the required smooth connection we are looking for,  we have obtained the same condition
which locates the discontinuity in $u$ of the profile for $\epsilon=0$.

Finally, computing the vector field $(G,H)=\left(w-g(u),v\right)$ at two given points $(u_\ell, v_\ast, w_\ast)\in \mathcal{W}_{-,0}^u$
and $(u_r, v_\ast, w_\ast)\in \mathcal{W}_{+,0}^s$ and taking a linear combination of them, we infer
\begin{align*}
	& \alpha(G,H) (u_\ell, v_\ast, w_\ast) +\beta (G,H)(u_r,v_\ast, w_\ast) \\
	& \ = ((\alpha + \beta) w_\ast + \alpha g(u_\ell) +\beta g(u_r), (\alpha + \beta) v_\ast ).
\end{align*}
Since $v_\ast<0$, see \eqref{vastnegative},  $u_r<u_\ell$, see \eqref{eq:ulurdiff}, and $g$ is monotone increasing, see \eqref{hypFG},
the two vectors are linearly independent which implies \eqref{nointer} and then, invoking Theorem \ref{thm:trans}, $N_{-,0}^u$ and
$N_{+,0}^s$ intersect transversally. 
Finally,  the heteroclinic orbit for the reduced system at $\epsilon=0$ persists for $\epsilon >0$ sufficiently small and the proof of
Theorem \ref{thm:main1} is complete.
Indeed, thanks to GSPT, we obtain the existence of the perturbed manifolds $N_{-,\epsilon}^u$ and $N_{+,\epsilon}^s$,
and their trasversality along a transversal heteroclinic orbit  is given by Theorem \ref{theo:a5}
(details on the construction of the manifolds $N_{-,\epsilon}^{u}$ and $N_{+,\epsilon}^{s}$ can be found in \cite{SZMOLYAN}).

\begin{example}[Hamer model]
The same conclusion holds for the viscous Hamer model, for which $R = \left[ \delta - \frac{1}{2} , 0\right] \times [a,b]$ and $d=1$.
Thus we refer to Theorem \ref{reduced} to check the transversality of $N_{+,0}^s$ and $N_{-,0}^u$.
As before, the condition of Theorem \ref{reduced} is equivalent to find $(v_0,w_0) \in R$ such that points $(\sqrt{1+2v_0},v_0,w_0)\in \mathcal{W}_{-,0}^u$
and $(-\sqrt{1+2v_0},v_0,w_0)\in \mathcal{W}_{+,0}^s$ are connected by 
\[
	u_0(x) = -u_\ell \tanh \left(\tfrac{1}{2}u_\ell x \right).
\]
As said before, the two maximal solution of the reduced systems $w_\pm$ are monotone decreasing, namely $w_+$ decreases from zero toward -1 at $+\infty$ and $w_-$ decreases from 1 at $-\infty$ toward 0. Hence,  we can conclude that the only point in $R$ verifying the conditions of Theorem \ref{reduced} is $(v_0,w_0)=(v^*,0)$. 
 The vector field $(G,H)=\left(w-u,v\right)$ computed at $(u_\ell,v^*,0) \in \mathcal{W}_{-,0}^u$ is given by $(-u_\ell,v^*)$, and at $(u_r,v^*,0) \in \mathcal{W}_{+,0}^s$
 is given by $(-u_r,v^*)=(u_\ell,v^*)$. Since $v^*\in \left[ \delta - \frac{1}{2} , 0\right]$ is not zero,    the two vectors are linearly independent and this gives \eqref{nointer}.
We stress that  $w=0$, guaranteeing the existence of the smooth connection between $u_\ell$ and $u_r$,
 is the same condition found in \cite{Kawashima} to locate the admissible sub--shock for the linear, inviscid Hamer model.
\end{example}

\section{Viscous radiating profile for small shocks}\label{sec:smallshocks}

In the previous sections we proved the existence of viscous radiating profiles for sufficiently small viscosity and large shocks,
in particular under the hypothesis that the profiles in the variable $u$ has a sub--shock.
 Here we complement such result by proving this existence in the case of $O(1)$-viscosity, fixed to be equal to 1, and sufficiently
 small shocks via a bifurcation argument with respect to the strength  of the shock.
This leads to the following system \eqref{genintro} with $\epsilon=1$, which is
\begin{equation}\label{PDE2}
	\left\{\begin{aligned}
	\partial_t u+ \partial_xf(u) - \partial_{xx}u & = \partial_x v \\
	v - \partial_{xx} v & = \partial_{x} g(u) 
	\end{aligned}\right.
\end{equation}
where we assume always hypotheses \eqref{hypFG}. 

We shall  prove the existence of solutions to \eqref{PDE2} in form of  travelling waves, i.e.\ $(u,v)(x,t)=\left(u(\xi),v(\xi)\right)$ 
where $\xi:=x-ct$ (with a slight abuse of notation), where the propagation speed is given by the Rankine--Hugoniot condition \eqref{RH}
and the Lax condition \eqref{Lax} is satisfied.
Since the flux $f$ is convex, the latter reduces to the inequality $u_+<u_-$.
The bifurcation parameter $\delta$ is defined as
\begin{equation*}\label{delta}
	\delta:= u_+ - u_- < 0,
\end{equation*}
and it is assumed to be sufficiently small.

In what follows, for the sake of simplicity, we shall focus on  the Hamer model, namely when 
$f(u) =\tfrac12\,u^2$ and $g(u) =u$, being the general case presented above analogous.
Note that, with such choices the Rankine-Hugoniot condition becomes
\begin{equation}\label{c}
		c=\frac{1}{2}\cdot\frac{u_+^2- u_-^2}{u_+ - u_-}=\frac{1}{2}(u_+ +u_-).
\end{equation}
and the Lax condition takes the form
\begin{equation}\label{lax3}
	u_+<c<u_-.
\end{equation}

\begin{theorem}\label{thm:main2}
Assume hypothesis \eqref{c}--\eqref{lax3}. 
Moreover, let the states $u_\pm$ be such that $|u_+-u_-|$ is sufficiently small.
Then the viscous Hamer system
\begin{equation*}\label{PDE22}
	\left\{ \begin{aligned}
	\partial_t u + \partial_{x}\left(\tfrac{1}{2} u^2 \right) - \partial_{xx} u & = \partial_{x}v \\
	v -\partial_{xx}v & = u,	
	\end{aligned} \right.
\end{equation*}
supports propagating fronts with speed $c$ given by equality \eqref{c}.
\end{theorem}

The remaining part of this section is devoted to the proof of this result. To this aim, let us start by recalling the dynamical system solved by the profile:
\begin{equation*}
	\left\{\begin{aligned}
	-cu' + \left(\tfrac12\,u^2\right)' - u'' & = v' \\
	v -v'' & = u'.
	\end{aligned}\right.
\end{equation*}
Integrating the first equation in $d \xi$, $\xi = x-ct$, from $\pm \infty$ we get:
\begin{equation*}
	\left\{\begin{aligned}
	u' & =\tfrac12\left( u^2-u_{\pm}\right)  -c(u-u_{\pm}) -v \\
	v -v'' & = u'.
	\end{aligned}\right.
\end{equation*}
Moreover, we can express the speed of the wave $c$ in \eqref{c} as a function of $u_+$ and $\delta$ as follows
\begin{equation*}
	c= \tfrac{1}{2}(2u_+-\delta)
\end{equation*}
and therefore
\begin{equation*}
	\left\{\begin{aligned}
	u' & = \tfrac{1}{2}(u-u_+)[(u-u_+)+\delta] -v \\
	v -v'' & = u'.
	\end{aligned}\right.
\end{equation*}
With the notation $\tilde{u}=u-u_+$ and $z= \tilde{u}+v'$, the dynamical system becomes:
\begin{equation}\label{ODE}
	X'=F(X;\delta),
\end{equation}
where  $X=(z,v,\tilde{u})$ and
\begin{equation*}
	F(z,v,\tilde{u};\delta)=\Bigl(v,z-\tilde{u}, \tfrac{1}{2}\tilde{u}^2 + \tfrac{1}{2}\tilde{u} \delta -v\Bigr).
\end{equation*}
Hence we observe that $p_1= (0,0,0)$ is a  critical point for any $\delta$, which corresponds to the point $(u_+,0,u_+)$ in the original variables. 
Moreover, depending on $\delta$, we have two different situations:
\begin{itemize}
\item[1)] if $\delta =0$ then $p_1$ is the only (trivial) critical point;
\item[2)] if $\delta < 0$, then a second critical point  $p_2=(-\delta,0, -\delta)$ bifurcates from the trivial one, the latter corresponds to the point $(u_-,0,u_-)$ in the original variables.
\end{itemize}

We want to transform \eqref{ODE} into its normal form performing a center manifold reduction, and prove a transcritical bifurcation
that occurs at $\delta=0$.
This will then imply the existence of the desired heteroclinic orbit, connecting $p_2$ at $- \infty$ to $p_1$ at $+\infty$. 
To this end we start by rewriting \eqref{ODE} with respect to the eigenbasis of the linearized system at the trivial critical point for $\delta=0$. 

Let $A=A(z,v,\tilde{u};\delta)$ be the Jacobian of $F$ at $(z,v,\tilde{u})$, namely 
\begin{equation*}
	A(z,v,\tilde{u};\delta)= \begin{pmatrix} 0 & 1 & 0 \\ 1 & 0 & -1 \\  0 & -1 & \tilde{u}+ \tfrac{1}{2}\delta \end{pmatrix}
\end{equation*}
The characteristic equation  at $p_1 = (0,0,0)$  is given by
\begin{equation*}
	-\lambda^3 + \tfrac{1}{2}\delta \lambda^2 + 2 \lambda - \tfrac{1}{2}\delta =0.
\end{equation*}
We study the real part of the eigenvalues at $p_1$ using the information coming from $\textrm{tr}\,(A|_{p_1})$ and $\det(A|_{p_1})$.
Since $\textrm{tr}\,(A|_{p_1})= \lambda_1+ \lambda_2+ \lambda_3 = \tfrac{1}{2}\delta < 0$ and $\det(A|_{p_1})= \lambda_1\lambda_2\lambda_3 = -\tfrac{1}{2}\delta > 0$,
we have two possibilities:
\begin{itemize}
\item[1)]  one real positive eigenvalue and two complex and conjugates  eigenvalues with negative real part;
\item[2)] two real negative  eigenvalues and one real positive eigenvalue.
\end{itemize}
Therefore at $p_1$ we have two stable direction and one unstable direction.
For $p_2 = (-\delta,0, -\delta)$ the situation is clearly reversed, being $\textrm{tr}\,(A|_{p_2})= -\tfrac{1}{2}\delta > 0 $ and $\det(A|_{p_2})= \tfrac{1}{2}\delta < 0$.
Hence at $p_2$ we have  two unstable direction and  one stable direction.

For $\delta =0$ the Jacobian at $p_1$ reduces to 
\begin{equation*}
	A(0,0,0;0)= \begin{pmatrix}
					0 & 1 & 0 \\ 
					1 & 0 & -1 \\ 
					0 & -1 & 0
			\end{pmatrix}
\end{equation*}
with eigenvalues $\lambda_1=0, \lambda_2 = \sqrt{2}, \lambda_3 =-\sqrt{2}$ and corresponding eigenvectors given
by $(1,0,1)$, $(-1,\sqrt{2},1)$ and $(1,\sqrt{2},-1)$.
Therefore, setting $Y:=(w_1,w_2,w_3)$, the desired change of basis is defined by the explicit matrix 
\begin{equation*}
	C =\begin{pmatrix} 1 & -1 & 1 \\ 0 & \sqrt{2} & \sqrt{2} \\  1 & 1 & -1 \end{pmatrix}
\end{equation*}
giving raise to
\begin{equation}\label{coo}
	Y=C^{-1}X.
\end{equation}
As a consequence, our original system \eqref{ODE} becomes
\begin{equation*}
	Y' =BY + C^{-1}\tilde{F}(CY;\delta),
\end{equation*}
where $B:=C^{-1}AC$ and $\tilde{F}(X; \delta) = F(X;\delta)- A(0,0,0) \cdot X$, namely:
\begin{align} \label{DEC}
	\frac{d}{dx}\begin{pmatrix} w_1 \\ w_2 \\ w_3 \end{pmatrix}
	& = \begin{pmatrix} 0 & 0 & 0 \\  0 & -\sqrt{2} & 0 \\ 0 & 0 & +\sqrt{2} \end{pmatrix}
		\begin{pmatrix} w_1 \\  w_2 \\  w_3  \end{pmatrix} \nonumber \\
	&\quad + \frac{1}{4}\begin{pmatrix}
		(w_1+w_2-w_3)^2 + \delta(w_1+w_2-w_3) \\ 
		\tfrac{1}{2}(w_1+w_2-w_3)^2 + \tfrac{1}{2}\delta(w_1+w_2-w_3) \\ 
		- \frac{1}{2}(w_1+w_2-w_3)^2 - \tfrac{1}{2}\delta(w_1+w_2-w_3)
		\end{pmatrix} 
\end{align}
In this way the system consists of a linear part and a perturbation.
The \emph{Center Manifold Theorem} 
guarantees the existence  for $\delta$ sufficiently small  of two $C^1$-functions
$\psi_2(w_1,\delta)$ and $\psi_3(w_1,\delta)$ such that 
\begin{equation*}
	w_2= \psi_2(w_1,\delta)\,,\quad w_3=\psi_3(w_1,\delta),
\end{equation*}
and the following tangency conditions hold:
\begin{equation}\label{eq:tang}
	\begin{aligned}
		&\psi_2(0,0)= \psi_3(0,0)= 0; \\
		&\frac{\partial \psi_2}{\partial w_1} (0,0) = \frac{\partial \psi_3}{\partial w_1} (0,0)= 0; \\
		&\frac{\partial \psi_2}{ \partial \delta}(0,0) = \frac{\partial \psi_3}{ \partial \delta}(0,0)=0.
	\end{aligned}
\end{equation}
Thus for $\delta < 0$ small enough, $\psi_2$ and $\psi_3$ are tangent in $(0,0)$ to the plane identified by $(w_1,\delta)=(0,0)$.
Moreover, the Center Manifold Theorem shows that, again for $\delta < 0$ small, in a neighborhood of the non hyperbolic critical
point $(0,0,0)$, the original system \eqref{ODE}, rewritten as in \eqref{DEC}, is $C^1$ topologically conjugate to the following decoupled system:
\begin{equation}\label{G}
	\left\{\begin{aligned}
	w_1'&= \tfrac{1}{4}\left\{w_1+ \psi_2(w_1,\delta) -\psi_3(w_1,\delta)\right\}^2\\
		&\quad + \tfrac{1}{4} \delta\left\{w_1+\psi_2(w_1,\delta)-\psi_3(w_1,\delta)\right\}\\
	w_2'&= -\sqrt{2}\,w_2, \\
	w_3'&=+ \sqrt{2}\,w_3.
	\end{aligned}\right. 
\end{equation}
Thus, to study the qualitative behaviour of the flow given by $\eqref{ODE}$ in a neighborhood of the non hyperbolic critical point $(0,0,0)$
and $\delta<0$ sufficiently small, we are reduced to the flow on the center manifold, which is given by $\eqref{G}_{1}$.
We want to apply \emph{Sotomayor Theorem} \cite[p.\ 338]{PERKO} to see which type of bifurcation occurs for $\delta=0$.
For completeness we report the statement below.

\begin{theorem}[Sotomayor]
Suppose that $F(X_{0}, \delta_{0})=0$ and that the matrix $A=JF(X_{0},\delta_{0})$ has a simple eigenvalue $\lambda=0$
with eigenvector $v$ and that $A^T$ has an eigenvector $w$ corresponding to the same eigenvalue.
Furthermore, suppose that $A$ has one eigenvalue with positive real part and one with negative real part and that the following conditions are satisfied:
\begin{equation}\label{w}
	\begin{aligned}
		&w^T \cdot  F_\delta(X_0,\delta_0) = 0;\\
		&w^T\cdot [DF_\delta(X_0,\delta_0)v] \neq 0;\\
		&w^T \cdot [D^2F(X_0,\delta_0)](v,v) \neq 0.
	\end{aligned}
\end{equation}
Then there is a smooth curve of equilibrium points of $X'=F(X;\delta)$ in $\mathbb{R}^3 \times \mathbb{R}$
through $(X_0,\delta_0)$ and tangent to $\mathbb{R}^3 \times \lbrace \delta_0 \rbrace$.
Depending on the sign in \eqref{w}, there are no equilibrium points near $X_0$ when $\delta < \delta_0$
(or when $\delta > \delta_0$) and there are two equilibrium points near $X_0$ when $\delta > \delta_0$ (or when $\delta< \delta_0$).
The two equilibria of $F(X;\delta)$ are hyperbolic and have stable manifolds of dimension one and two, respectively;
i.e., the system $X'=F(X; \delta)$ experiences a 
transcritical bifurcation at $X_0$ as the parameter $\delta$ passes through $\delta_0$.
\end{theorem}

Introducing the same notation of the theorem, with $\hat F(w_1,w_2,w_3;\delta)$ denoting the vector field
\begin{equation*}
	\begin{pmatrix}
	\frac{1}{4}\left\{w_1+ \psi_2(w_1,\delta) -\psi_3(w_1,\delta)\right\}^2 + \frac{1}{4} \delta\left\{w_1+\psi_2(w_1,\delta)-\psi_3(w_1,\delta) \right\} \\ 
	-\sqrt{2}\,w_2 \\ 
	+\sqrt{2}\,w_3
	\end{pmatrix}
\end{equation*}
we clearly have 
\begin{equation*}
	\hat F(0,0,0;0)= \begin{pmatrix} 0 \\  0 \\  0 \end{pmatrix} \quad\textrm{and}\quad
	J\hat F(0,0,0;0)= B =\sqrt{2} \begin{pmatrix} 0 & 0 & 0 \\  0 & -1 & 0 \\  0 & 0 & +1 \end{pmatrix}.
\end{equation*}
The eigenvector associated to $\lambda=0$ is $v=(1,0,0)$ and since $B=B^T$ we have that the corresponding eigenvector
to $\lambda=0$ for this matrix is the same, namely $w^T=(1,0,0)$. 
Moreover, in view of the  tangency conditions \eqref{eq:tang}, a direct calculation shows
\begin{equation*}
	\hat F_{\delta}(0,0,0;0) = (0,0,0);\qquad
	D\hat F_{\delta}(0,0,0;0) =\frac14 \begin{pmatrix}
		1 & 0 & 0 \\  0 & 0 & 0 \\  0 & 0 & 0
	\end{pmatrix}.
\end{equation*}
Therefore $w^T \cdot \hat F_{\delta}(0,0,0;0) =0$ and $w^T \cdot D\hat F_{\delta}(0,0,0;0)v = 1/4 \neq 0$. 
Finally we have to compute $w^T \cdot[D^2\hat F(0,0,0;0)(v,v)]$ and we have to check that is different from zero.
The computation of $[D^2\hat F(0,0,0;0)(v,v)]$ gives  the vector $(1/2,0,0)$ and thus $w^T \cdot[D^2\hat F(0,0,0;0)(v,v)] \neq 0$. 
In view of  \emph{Sotomayor Theorem}, these three conditions imply that the original system $\eqref{ODE}$ experiences
a transcritical bifurcation at the equilibrium point $(0,0,0)$ as the parameter $\delta$ varies through the bifurcation value $\delta=0$. 

The presence of the parameter $\delta$ is only in the first component $\hat F^1$ of the vector field $F$, for which,
evaluating $\hat F^1_{w_1w_1}(0,0,0;0) = \tfrac12$ and  $F^1_{\delta w_1}(0,0,0;0) =\tfrac14$,
we obtain the normal form (among others, see \cite[Theorem 1.3]{an introduction}): 
\begin{equation*}
	w_1'= \tfrac{1}{4}\delta w_1 + \tfrac{1}{4} w_1^2 = \tfrac{1}{4} w_1(w_1 + \delta) .
\end{equation*}
Since $\delta < 0$, the trivial equilibrium $w_1=0$ is stable and $w_1=-\delta > 0$ is unstable, and the same is true for  the one--dimensional  center manifold reduction given by 
\begin{equation}\label{eq:G1}
	\begin{aligned}
	w_1'=&\tfrac{1}{4}\left\{w_1+ \psi_2(w_1,\delta) -\psi_3(w_1,\delta)\right\}^2 \\
		&+ \tfrac{1}{4} \delta\left\{w_1+\psi_2(w_1,\delta)-\psi_3(w_1,\delta)\right\}.
	\end{aligned}
\end{equation}
As a consequence, there exists an heteroclinic conncetion between $w_1=-\delta$ at $-\infty$ and $w_1=0$ at $+\infty$ for that equation. 
Then \eqref{ODE} is locally topologically equivalent to \eqref{eq:G1} augmented with the two  linear equations 
\begin{equation*}
	w_2'= -\sqrt{2}\,w_2,\  w_3'=+\sqrt{2}\,w_3,
\end{equation*}
 namely, system \eqref{G} (see \cite[Theorem 5.4, p.\ 159]{K}). 
The number of positive and negative eigenvalues is preserved, as well as the  trajectories of the two dynamical systems. 
We have two unstable and one stable directions for the equilibrium point $(-\delta,0,0)$ and two stable directions and only one unstable at $(0,0,0)$.
The two--dimensional unstable eigenspace of $(-\delta,0,0)$ and the two--dimensional stable eigenspace of $(0,0,0)$ intersect along the direction
given by the vector $(1,0,0)$, namely the tangential direction of the center manifold in a neighborhood of the trivial equilibrium point.
This implies the heteroclinic actually exists for \eqref{G}.
Using \eqref{coo} to go back to the original variables, we have:
\begin{itemize}
\item[1)] the point $(0,0,0)$ is mapped to $(u_+,0,u_+)$,
\item[2)] the point $(-\delta,0,0)$ is mapped to $(u_-,0,u_-),$ 
\end{itemize}
and we have finally proved the existence of an heteroclinic orbit between these two points as stated in Theorem \ref{thm:main2},  thanks to the aforementioned local topological equivalence. 

\appendix
\section{Geometric Singular Perturbation Theory}\label{app:GSPT}

\subsection{Invariant Manifold Theorems}

In this section we recall the main results about \emph{Geometric Singular Perturbation Theory} developed in \cite{Fenichel, SZMOLYAN}. Let us consider the vector field $X_\epsilon \times \lbrace 0 \rbrace$ defined in the following way:
\begin{equation*}
\left\{
\begin{split}
&x' = \epsilon f(x,y,\epsilon) \\
&y' = g(x,y,\epsilon) \\
&\epsilon' = 0 
\end{split} \right.
\end{equation*}
where $x \in \mathbb{R}^n$ are the slow variables, $y \in \mathbb{R}^k$ are the fast variables. The algebraic costraint $g(x,y,0) = 0$ defines the slow manifold $\mathcal{S}$ and correspondingly the two invertible branches $\mathcal{S}^{\pm}$. The following invariant manifold theorem describes the flow induced by $X_\epsilon \times \lbrace 0 \rbrace$ near $\mathcal{S}^\pm \times \lbrace0\rbrace$ for small $\epsilon$.

\begin{theorem}\label{existence}
Let $M$ be a $C^{r+1}$ manifold, $1 \leq r < \infty$.
Let $X_{\epsilon}$, $\epsilon \in (-\epsilon_0,\epsilon_0)$ be a $C^r$ family of vector fields on $M$, and let $\mathcal{S}$ be
a $C^r$ submanifold of $M$ consisting entirely of equilibrium points of $X_0$. Let $k^s$, $k^c$ and $k^u$ be fixed integers,
and let $K \subset \mathcal{S}^\pm$ be a compact subset such that $QX_0(m)$ has $k_s$ eigenvalues in the left half plane, $k_c$ eigenvalues on the
imaginary axis, and $k_u$ in the right half plane, for all $m \in K$.
Then:
\begin{itemize}
\item[1)] There is a $C^r$ center-stable manifold $\mathcal{C}^s$ for $X_{\epsilon} \times 0$ near $K$.
There is a $C^r$ center-unstable manifold $\mathcal{C}^u$ for $X_{\epsilon} \times 0$ near $K$.
There is a $C^r$ center manifold $C$ for $X_{\epsilon} \times 0$ near $K$.
\item[2)] There is a $C^{r - 1}$ family $\mathcal{F}^s = \lbrace \Psi^s(p) : p \in \mathcal{C}^s \rbrace $ of $C^r$ stable manifolds for $\mathcal{C}^s$ near $K$.
If $p \in M \times \lbrace \epsilon \rbrace$, then $\Psi^s(p) \in M \times \lbrace \epsilon \rbrace$.
Each manifold $\Psi^s(p)$ intersects $\mathcal{C}$ transversally, in exactly one point.
There is a $C^{r-1}$ family $\mathcal{F}^u = \lbrace \Psi^u(p): p \in \mathcal{C}^u \rbrace$ of $C^r$ unstable manifolds for $\mathcal{C}^u$ near $K$.
If $p \in M \times \lbrace\epsilon\rbrace$, then $\Psi^u(p) \in M \times \lbrace\epsilon\rbrace$.
Each manifold $\Psi^u(p)$ intersects $\mathcal{C}$ transversally, in exactly one point.
\item[3)] Let $K_s < 0$ be larger then the real parts of the eigenvalues of $QX_0(m)$ in the left half plane, for all $m\in K$.
Then, there is a constant $C_s$ such that if $p \in \mathcal{C}^s$ and $q \in \Psi^s(p)$, then
\begin{equation*}
	d(p \cdot \bar{x}, q \cdot \bar{x}) \leq C_s e^{K_s\bar{x}}d(p,q)
\end{equation*}
for all $\bar{x} \geq 0$ such that $p \cdot [0,\bar{x}] \subset \mathcal{C}^s$.
Let $K_u > 0$ be smaller than the real parts of the eigenvalues of $QX_0(m)$ in the right half plane, for all $m \in K$.
Then there exists a constant $C_u$ such that if $p \in \mathcal{C}^u$ and $q \in \Psi^u(p)$, then
\begin{equation*}
	d(p \cdot \bar{x}, q \cdot \bar{x}) \leq C_u e^{K_u \bar{x}}d(p,q)
\end{equation*}
for all $\bar{x} \leq 0$ such that $p \cdot [\bar{x}, 0] \subset \mathcal{C}^u$.
\item[4)] Let $S_H \subset S^\pm$ such that $QX_0(m)$ has not eigenvalue with zero real part. If $K \subset \mathcal{S}_H$, define for $(m, \epsilon) \in \mathcal{C}$,
\begin{equation*}
	\mathbf{F}_R(m) := \pi^S \left( \frac{\partial}{\partial \epsilon} \right) X_{\epsilon}(m)|_{\epsilon = 0}
\end{equation*}
and
\begin{equation*}
	X_{\mathcal{C}}(m, \epsilon) := \left\{\begin{aligned}
		& \epsilon^{-1}X_{\epsilon}(m) \times \lbrace 0 \rbrace, \; if \; \epsilon \neq 0\\
		& \mathbf{F}_R(m) \times \lbrace 0 \rbrace, \; if \; \epsilon =0
		\end{aligned}\right.
\end{equation*}
Then $X_{\mathcal{C}}$ is a $C^{r-1}$ vector field on $\mathcal{C}$ near $K \times \lbrace0\rbrace$.
\end{itemize}
\end{theorem}

\subsection{Reduced System}

The following lemma gives the reduced system in Theorem \ref{existence} of the previous section in local coordinates in which
$\mathcal{S}$ appears as graph of a function (see \cite{Fenichel}).
From now on, the symbols $D_i$, $i=1,2,3$ denote the derivative with respect to slow variables $D_x$, fast variables $D_y$ and $D_{\epsilon}$ respectively.

\begin{lemma}\label{lemma3.1}
Consider the system
\begin{equation*}
		\dot{x} = f(x,y,\epsilon),\qquad
		\dot{y} = g(x,y,\epsilon),
\end{equation*}
defined for $(x,y)$ in an open subset of $\mathbb{R}^n \times \mathbb{R}^k$, for $\epsilon$ near zero.
Let $y=h(x)$ be a function such that 
\begin{equation*}
	f(x,h(x),0) \equiv 0,
	\quad\textrm{and}\quad
	g(x,h(x),0) \equiv 0.
\end{equation*}
Suppose $(x_0,h(x_0))\in \mathcal{S_H}$, so that the matrix 
\begin{equation*}
	\begin{pmatrix} \alpha & \beta \\  \gamma & \delta \end{pmatrix}
	= \begin{pmatrix} 	\partial_x f(x_0, h(x_0),0)	& \partial_y f(x_0, h(x_0),0) \\ 
					\partial_x g(x_0, h(x_0),0)	& \partial_y g(x_0, h(x_0),0)  \end{pmatrix}
\end{equation*}
has rank $k$. Let $\nu = \partial_x h(x_0)$.
Then the projection $\pi^\mathcal{S}$ = $\pi^\mathcal{S}(x_0,u(x_0))$ is multiplication by the matrix
\begin{equation*}
	\mathbb{A}=\begin{pmatrix} I+ \beta(\delta- \nu \beta)^{-1}\nu 			& -\beta(\delta- \nu \beta)^{-1} 		\\ 
					\nu + \nu \beta(\delta- \nu \beta)^{-1}\nu	& -\nu \beta(\delta- \nu \beta)^{-1}	\end{pmatrix}
\end{equation*}
and the reduced system is given by:
\begin{equation*}
	\begin{pmatrix} \dot x \\ \dot y \end{pmatrix}
	= \mathbb{A} \begin{pmatrix} \partial_\epsilon f(x,h(x),0) \\  \partial_\epsilon g(x,h(x),0) \end{pmatrix}
\end{equation*}
\end{lemma}

\subsection{Local theory near equilibria of the reduced vector field}

The next theorem explains how the critical points of the reduced vector field at $\epsilon=0$ -- denoted by $\mathbf{F}_R$ --
are related to the critical points of the original vector field $X_{\epsilon}$ in the regime of small $\epsilon > 0$.

Given $\epsilon \in (-\epsilon_0, \epsilon_0)$, let $X_\epsilon$ be a $C^r$ family of vector fields on $\mathbb{R}^3$
and let $\mathcal{S}$ be a $C^r$ submanifold of $\mathbb{R}^3$ consisting of equilibrium points of $X_0$.
Moreover, let $\mathbf{u}\in \mathcal{S}_H$ be an equilibrium point of the reduced vector field $\mathbf{F}_R$.

\begin{theorem}[Theorem 12.1, \cite{Fenichel}]\label{teo3.3}
If $\lambda=1$ is not an eigenvalue of $T\mathbf{F}_R(m)$, then there exists $\epsilon_1 > 0$ and a $C^{r-1}$ family
of points $\mathbf{u}_\epsilon$ with $\epsilon \in (-\epsilon_1, \epsilon_1)$ such that $\mathbf{u}^0= \mathbf{u}$ and $\mathbf{u}_\epsilon$ is an
equilibrium point of $X_\epsilon$.
\end{theorem}

Let $\cdot t$ denotes a flow on the manifold $M$.
Given $V\subset M$, set
\begin{equation*}
	\begin{aligned}
	A^+(V)&:=\left\{ \mathbf{u} \in V\,:\, \overline{\mathbf{u} \cdot [0,+\infty)}\subset V\right\}\,,\\
        	A^-(V)&:=\left\{ \mathbf{u} \in V\,:\, \overline{\mathbf{u} \cdot (- \infty,0]}\subset V \right\}\,,\\
	I(V)&:= \left\{ \mathbf{u} \in V\,:\, \overline{\mathbf{u} \cdot (- \infty, \infty)}\subset V \right\}\,.
	\end{aligned}
\end{equation*}
The following result characterizes the local stable/unstable manifolds of normally hyperbolic invariant manifold of $\mathbf{F}_R$.

\begin{theorem}[Theorem 12.2, \cite{Fenichel}, Theorem 2.2, \cite{SZMOLYAN}]\label{s-u-mxr}
Under the hypothesis of Theorem \ref{teo3.3}, suppose $T\mathbf{F}_R(\mathbf{u}_0)$ has $j^u$ eigenvalues in the right half plane,
no eigenvalues on the imaginary axis and $j^s$ eigenvalues in the left half plane.
Suppose $Q\mathbf{F}_0(\mathbf{u}_0)$ has $k^u$ eigenvalues in the right half plane and $k^s$ eigenvalues in the left half plane.

Then, there exists $\epsilon_1>0$ such that
\begin{itemize}

\item[i.] there is a $C^{r-1}$-family of hyperbolic equilibrium points of $\mathbf{F}_\epsilon$, denoted by 
$\left\{\mathbf{u}_\epsilon: \epsilon \in (-\epsilon_1,\epsilon_1)\right\}$, such that $\lim\limits_{\epsilon\to 0}\mathbf{u}_\epsilon=\mathbf{u}_0$
and there is a family of neighborhoods of $\mathbf{u}_0$, denoted by $\left\{ U_\epsilon\,:\, \epsilon \in (-\epsilon_1,\epsilon_1) \right\}$,
such that $I(U_\epsilon) = \left\{ \mathbf{u}_\epsilon \right\}$ for any $\epsilon \neq 0$;

\item[ii.] there is a $C^{r-1}$-families of $(j^u+k^u)$-dimensional manifolds $\lbrace W_\epsilon^u : \epsilon \in (-\epsilon_1,\epsilon_1) \rbrace$
and $(j^s +k^s)$-dimensional manifolds $\left\{ \mathcal{W}_\epsilon^s : \epsilon \in (-\epsilon_1,\epsilon_1) \right\}$ such that
\begin{equation*}
	A^-(U_\epsilon)=W_\epsilon^u \quad\textrm{and}\quad  A^+(U_\epsilon) = \mathcal{W}_\epsilon^s \quad\forall\,\epsilon>0;
\end{equation*}

\item[iii.] the local stable and unstable manifolds of $\mathbf{u}_{\epsilon}$ for $\epsilon>0$ are given by
\begin{equation*}
	N_\epsilon^s = \bigcup_{p \in \mathcal{W}_\epsilon^s} \Psi_\epsilon^s(p)
	\quad\textrm{and}\quad
	N_\epsilon^u = \bigcup_{p \in \mathcal{W}_\epsilon^u} \Psi_\epsilon^u(p),
\end{equation*}
where $\lbrace \Psi_\epsilon^s(\mathbf{u}): \epsilon \in (-\epsilon_1,\epsilon_1), \mathbf{u} \in \mathcal{W}_\epsilon^s \rbrace$ is a $C^{r-1}$-family of
$k_1$-dimensional manifolds such that $\lbrace \Psi_\epsilon^s(m): m \in \mathcal{W}_\epsilon^s \rbrace$ for each $\epsilon>0$
is a positively invariant family of manifolds (the same for $\Psi_\epsilon^u(m)$).
\end{itemize}
\end{theorem}

\begin{theorem}[Theorem 3.1, \cite{SZMOLYAN}]
\label{theo:a5}
Let the manifolds $N_1$ and $N_2$ satisfy the assumptions of Theorem \ref{s-u-mxr}.
Assume that the manifolds
\begin{equation*}
	N_1^u:=\bigcup_{p\in \mathcal{W}_1^u} \Psi^u(p),\quad
	N_2^s:=\bigcup_{p\in \mathcal{W}_2^s} \Psi^s(p)
\end{equation*}
intersect transversally along the singular heteroclinic orbit.

Then, there exists $\epsilon_1>0$ such that there exists a transversal heteroclinic orbit of the singularly perturbed system 
\begin{equation*}
		\dot{x} = f(x,y,\epsilon),\qquad
		\epsilon \dot{y} = g(x,y,\epsilon),
\end{equation*}
connecting the manifolds $N_{1,\epsilon}$ and $N_{2,\epsilon}$ for $0<\epsilon<\epsilon_1$.
\end{theorem}

Finally, the last needed result concerns sufficient (transversality) conditions needed to obtain the connections
between the two branches of the involved invariant manifold.

\begin{theorem}[Theorem 4.1, \cite{SZMOLYAN}]\label{teo3.12}
Let the manifolds $N_1$ and $N_2$ satisfy the assumption of Theorem $\ref{s-u-mxr}$. Let $\phi(W_1^u)$ and $\phi(W_2^s)$
denote the $x$-coordinates of the manifolds $\mathcal{W}_1^u$ and $\mathcal{W}_2^s$, respecitvely.
Then the manifolds $N_1^u$ and $N_2^s$ intersect transversally in the points of the heteroclinic orbit $(x_0,y_0(x))$
if and only if there exist exactly $d-1$ linearly independent solutions $\xi \in T_{x_0}\phi(\mathcal{W}_1^u) \bigcap T_{x_0}\phi(\mathcal{W} _2^s)$
of the equation
\begin{equation*}
	(M,\xi)= 0\quad\textrm{where}\quad 
	M:= \int_{\mathbb{R}} \psi(\xi) \cdot \partial_{x}g\bigl(x_0, y_0(\xi)\bigr)\,d\xi.
\end{equation*}
The function $\psi$ is the unique (up to a scalar multiple) bounded solution of the adjoint equation 
\begin{equation*}
	\psi'= -\partial_y g(x_0,y_0(x))^T \psi.
\end{equation*}
\end{theorem}

\end{document}